\documentclass[12pt]{article}
\usepackage{amsmath,amssymb,amscd,latexsym,mathrsfs}
\usepackage[all]{xy}

\def\Hom{\mathrm{Hom}}
\def\Lan{\mathrm{Lan}}
\def\Set{\mathrm{Set}}
\def\Ab{\mathrm{Ab}}
\def\Ob{\mathrm{Ob}}
\def\Mor{\mathrm{Mor}}
\def\Comp{\mathrm{Ch}}

\def\mC{\mathscr C}
\def\mD{\mathscr D}
\def\mA{\mathcal A}
\def\II{\mathbb I}

\def\NN{\mathbb N}
\def\ZZ{\mathbb Z}
\def\coLim{\underrightarrow{\lim}}

\newtheorem{theorem}{Theorem}
\newtheorem{lemma}{Lemma}[section]
\newtheorem{propos}{Proposition}
\newtheorem{corollary}{Corollary}
\newtheorem{definition}{Definition}

\newtheorem{example}{Example}

\def\leq{\leqslant}
\def\geq{\geqslant}

\begin{document}

\begin{center}
{\Large Homology groups of cubical sets}
\\
\medskip
Ahmet A. Husainov\\
\end{center}

\date{}

\begin{abstract}
The paper is devoted to homology groups of cubical sets
with coefficients in contravariant systems of Abelian groups.
The study is based on the proof of the assertion that the homology groups
of the category of cubes with coefficients in the diagram of Abelian groups
are isomorphic to the homology groups of normalized complex of the cubical Abelian
group corresponding to this diagram.
The main result
shows that the homology groups of a cubical set with
coefficients in a contravariant system of Abelian groups are isomorphic to
the values of left derived functors of the colimit functor on this contravariant system.
This is used to obtain the isomorphism criterion for homology
groups of cubical sets with coefficients in contravariant systems, and also to
construct spectral sequences for the covering of a cubical set and for a morphism
between cubical sets.
\end{abstract}




\section{Introduction}\label{s1}

In this paper we study the homology groups of cubical sets with coefficients in contravariant systems of Abelian groups. This theme grew out of the classical works of J.P. Serre \cite{ser1958}, S. Eilenberg and S. MacLane \cite{eil1953}, where the homology of the cubical set of singular cubes of topological space
 was considered. Serre introduced homology groups with coefficients in local systems of Abelian groups.
It is been quite a long time.
Further development is related to applications. Homological methods were applied in the theory of image analysis \cite{nie2002}, which stimulated the development of the homology theory of cubical subsets of Euclidean spaces \cite{kac2003} and methods for computing these groups  \cite{har2014} - \cite{pil2015}. Homology and homotopy of cubical sets are used in the noncommutative topology \cite{bro2011}.

Starting with the thesis of E. Goubault \cite{gou1995} devoted to the homology of higher-dimensional automata, papers began to appear containing applications of the homology theory of cubical and precubical sets for studying mathematical models for computational processes and systems.

It became clear that for the needs of the theory of computational processes and systems it is not enough to consider homology groups with coefficients in local systems.

Similar problems arose in the directed topology \cite[Problem (M4)]{faj2006}.
In the monograph of M. Grandis \cite{gra2009}, a homology theory of cubical sets with coefficients in ordered Abelian groups was developed and homology groups of directed topological spaces were introduced.
We studied homology groups of asynchronous systems and Petri nets with coefficients in modules on free partially commutative monoids \cite{X2004}. It was proved \cite{X20082} - \cite{X2013} that these groups can be calculated as homology groups of precubical sets with coefficients in some contravariant systems.

In \cite{X20132}, the adjoint functors between the category of cubical sets and the category of asynchronous systems were constructed. This indicates that it is time to study the homology of cubical sets with coefficients in contravariant systems that may not be local.
Here by local systems we mean contravariant systems consisting of isomorphisms.

We consider an arbitrary abstract cubical set $X$ with a functor from the category of singular cubes of the cubical set $X$ to the category of Abelian groups. This functor is called the contravariant system on $X$. Our goal is to study the homology groups of a cubical set with coefficients in contravariant systems.

The main problem can be formulated as follows.
It is known \cite[Application 2]{gab1967} that the homology groups $H_n(X,F)$
of a simplicial set $X$ with coefficients in the contravariant system  $F$ are isomorphic
to the values $\coLim^{(\Delta/X)^{op}}_n F$ of the left derived functors of the colimit functor.
The same is true for contravariant systems of Abelian groups on semisimplicial sets (see, for example \cite[Proposition 1.4]{X1997} where the dual assertion is proved).
In \cite[Theorem 4.3]{X2008}, this assertion was proved for precubical sets.

The homology groups of a cubical set, in contrast to the homology groups of simplicial set, are defined by means of a normalized complex.
 This makes it difficult to study them. Will the homology groups of a cubical set be equal to the values of the derived functors of the colimit functor?

We give a positive answer.
First, we construct a projective resolution in the category of cocubical Abelian groups whose tensor product on an arbitrary cubical Abelian group is isomorphic to the normalized complex of this cubical Abelian group. On the basis of this, we establish that the homology groups of the category of cubes with coefficients in a cubical Abelian group are isomorphic to the homology groups of the normalized chain complex of this cubical Abelian group.
Similar results were known for homology groups
of simplicial Abelian groups, semisimplicial Abelian groups, precubical groups, and
cyclic objects in the category of Abelian groups \cite[Example 2.1]{X2002}.

Then we prove the main theorem of this paper that the homology groups
$H_n (X, F)$ of a cubical set with coefficients in the contravariant system of Abelian groups
$F: (\Box / X)^{op}\to \Ab$ are isomorphic to the values $ \coLim^{(\Box/X)^{op} }_n F$ of the left derived of the colimit functor.

We apply the main theorem to obtain an isomorphism criterion for the homology groups of cubical sets, and also to construct spectral sequences for the covering of a cubical set and for a morphism of cubical sets.

\section{Preliminaries}

Let $\Set$ be the category of sets and maps and let $\Ab$ be the category
of Abelian groups and homomorphisms.
We denote by $\ZZ(-): \Set \rightarrow \Ab$ the functor
 which assigns to each set $E$ the free Abelian group $\ZZ(E)$ with basis $ E $ and to each map
 $f: E_1 \rightarrow E_2 $ the homomorphism $ \ZZ (f): \ZZ (E_1) \rightarrow \ZZ (E_2) $, which extends this mapping.
Denote by $\II= \{0, 1\}$ the set ordered by the relation
$0< 1$. Let $\ZZ$ be the set or additive group of integers and
$\NN$ be the set of nonnegative integers.

For any category $\mA$, we denote by $\mA^{op}$ the category opposite to $\mA$.
 For arbitrary objects $a,b\in \mA$ we denote by
$\mA(a,b)$ the set of morphisms $a\rightarrow b$.
For any morphisms $\alpha: a'\to a$, $\beta: b\to b'$ there are maps
$\mA(\alpha,\beta): \mA(a,b)\to \mA(a',b')$ defined as
$\mA(\alpha,\beta)(\gamma)= \beta\gamma\alpha$.
If $\mD$ is a small category, then functors
 $\mD\to \mA$ are called  {\it diagrams} of objects on $\mD$ in $\mA$.
We denote by $\mA^{\mD}$ the category
of diagrams on $\mD$ in $\mA$ and natural transformations between them.

\subsection{Tensor product of diagrams}

For any small category $\mD$ and cocomplete additive category $\mA$ there is a bifunctor of tensor product
$$
\otimes: \Ab^{\mD}\times \Ab^{\mD^{op}} \to \Ab,
$$
which can be characterized by means of isomorphism
$$
\Ab(G\otimes F, A)\stackrel{\cong}\to \Ab^{\mD}(G, \Hom(F(-),A))
$$
that must be natural in each argument.
Here $\Hom(F(-),A): \mD\to \Ab$ is the functor with values
 $\Ab(F(c),A)$ for all $c\in \Ob\mD$. Its values on the morphisms
$\alpha: a\to b$ are equal to maps
$\Ab(F(\alpha),1_A): \Ab(F(b),A)\to \Ab(F(a),A)$.

For an arbitrary $c\in \Ob\mD$, denote by $\ZZ h^c: \mD\to \Ab$ the
composition of functors $\mD \stackrel{h^c}\to \Set\stackrel{\ZZ(-)}\to \Ab$
where $h^c=\mD(c,-)$.
It is well known (see, for example, \cite[Lemma 3.2]{X2008}) that
there is a natural isomorphism $\xi_c: \ZZ h^c\otimes F\stackrel{\cong}\to F(c)$.
The natural property means that for every morphism $ \alpha: a \to b $ the following diagram is commutative, shown in Fig. 1.
\begin{figure}[h]
$$
\xymatrix
{
\ZZ h^a\otimes F \ar[rr]^{\xi_a} && F(a)\\
\ZZ h^b\otimes F \ar[u]^{\ZZ h^{\alpha}\otimes 1_F} \ar[rr]_{\xi_b}
 && F(b)\ar[u]_{F(\alpha)}
}
$$
\centerline{\small Fig. 1. Natural isomorphism of functors $\ZZ h^{(-)}\otimes F \to F$.}
\end{figure}
For fixed  $F\in \Ab^{\mD^{op}}$, the functor
 $(-)\otimes F$ commutes with the colimits. For each
 $G\in \Ab^{\mD}$ the functor  $G\otimes (-)$ commutes with the colimits.

\subsection{Derived functors of the colimit}
Let $\mD$ be a small category.
Denote by $\Delta_{\mD}\ZZ: \mD\to \Ab$ the diagram of Abelian groups taking the
constant values $\ZZ$ on objects, and the values $\Delta_{\mD}\ZZ(\alpha)=1_{\ZZ}$ on the morphisms $\alpha\in \Mor\mD$ where
$1_{\ZZ}: \ZZ\to \ZZ$ is the identity homomorphism of the additive group of integers.

The category $\Ab^{\mD^{op}}$ is Abelian. It is cocomplete and has enough projectives.
The colimit functor $\coLim^{\mD^{op}}: \Ab^{\mD^{op}}\to \Ab$ is right exact.
Hence, the colimit functor has left derived functors
$\coLim^{\mD^{op}}_n: \Ab^{\mD^{op}}\to \Ab$. The description of the chain complex whose homology groups are naturally isomorphic to $\coLim^{\mD^{op}}_n F$ is given in
\cite[Application 2]{gab1967} and in \cite[Definition 3.1]{X2008}.
For the definition of the derived functors of the colimit, we can take the following assertion:
\begin{propos}\label{projcolim}\cite{obe1967}
For each projective resolution $P_*\to \Delta\ZZ$ in the category $\Ab^{\mD}$, there are
isomorphisms $\coLim_n^{\mD^{op}}F\cong H_n(P_*\otimes F)$,
which are natural in $F\in \Ab^{\mD^{op}}$.
\end{propos}

\section{Homology of the category of cubes}

In this section, we recall cubes and cubical objects, and
construct a projective resolution of the cocubical Abelian group $ \Delta _{\Box} \ZZ$.
We prove that the tensor product of this resolution on an arbitrary cubical Abelian group $F$ is isomorphic to the normalized complex of  $F$.
It follows that the homology groups of this normalized complex are isomorphic
to the homology groups $\coLim^{\Box^{op}}_n F$ of the category $\Box^{op}$
with coefficients in $F$.

\subsection{The category of cubes}
For an arbitrary $n\in \NN$ we consider
the partially ordered set $\II^n=\{0,1\}^n$ equal to the $n$th Cartesian power of the
 totally ordered set $\II=\{0,1\}$. For $n=0$ the set $\II^0$ consists of the single element
 $\emptyset$.
The partially ordered set $\II^n$ is called the
$n$-dimensional {\it cube}.

The objects of the {\it category of cubes} $\Box$ are the cubes
$\II^0$, $\II^1$, $\II^2$,  $\cdots$.

Morphisms of the category $\Box$ are nondecreasing maps $\II^p\to \II^q$
that can be got by composition of the {\it face morphisms}
$\delta_i^{k,\tau}: \II^{k-1}\rightarrow \II^k$ and of {\it degeneracy morphisms}
$\epsilon_i^k: \II^k\to \II^{k-1}$ defined for $k\geq 1$, $1\leq i\leq k$,
$\tau\in \{0,1\}$, respectively by
\begin{gather}
\delta_i^{k,\tau}(x_1, \dots, x_{k-1})=
(x_1, \dots, x_{i-1}, \tau, x_i, \dots, x_{k-1}),\\
\epsilon_i^k(x_1, \cdots, x_k)=(x_1, \cdots, x_{i-1}, x_{i+1}, \cdots, x_k).
\end{gather}
In particular, $\delta_1^{1,0}(\emptyset)=0$, $\delta_1^{1,1}(\emptyset)=1$,
$\epsilon^1_1(x_1)=\emptyset$, for all $x_1\in \II$.

Sufficiently complete information on the category of cubes can be found in \cite{bro2011}, \cite{gra2003}, and \cite{jar2006}. In particular \cite{bro2011}, its objects can be considered as
Euclidean cubes $[0,1]^n$.
We have the equalities
\begin{gather}
\label{rel1}\delta_j^{n,\beta}\delta_i^{n-1,\alpha}=
\delta_i^{n,\alpha}\delta_{j-1}^{n-1,\beta}\quad  (1\leq i< j \leq n,
\alpha\in \II, \beta\in \II);\\
\label{rel2}
\epsilon^{n-1}_j \epsilon^n_i= \epsilon^{n-1}_i\epsilon^n_{j+1}\quad (1\leq i\leq j\leq n-1,
n\geq 2);\\
\label{rel3}
\epsilon^{n+1}_j\delta^{n+1,\alpha}_i=\left\{
\begin{array}{lc}
\delta^{n,\alpha}_i\epsilon^{n}_{j-1}, & (1\leq i<j\leq n+1, \alpha\in \{0,1\}), \\
\delta^{n,\alpha}_{i-1}\epsilon^{n}_j, & (1\leq j<i\leq n+1, \alpha\in \{0,1\}),\\
1_{\II^n}, & i=j.
\end{array}
\right.
\end{gather}

By \cite[Lemma 4.1]{gra2003}, every morphism
$f: \II^k\to \II^n$ in $\Box$ has a unique decomposition of the form
\begin{equation}\label{grandis}
f= \delta^{n,\tau_1}_{j_1}\cdots \delta^{k-r+1,\tau_s}_{j_s}
\epsilon^{k-r+1}_{i_1}\cdots \epsilon^k_{i_r},
\quad
\begin{array}{l}
 1\leq i_1< \cdots <i_r\leq k,\\
 n\geq j_1>\cdots >j_s\geq 1,\\
 k-r=n-s\geq 0.
\end{array}
\end{equation}

A morphism in $\Box$ is {\it injective} if it is given by an injection.
The category $\Box$ contains the subcategory $\Box_+$ such that
 $\Ob(\Box_+)= \Ob(\Box)$ and morphisms in $\Box_+$ are injective morphisms of
$\Box$. This subcategory is generated by face morphisms $\delta^{k,\tau}_i$
and can be given by the relations (\ref{rel1}).

\subsection{Cubical sets}
Let $\mA$ be a category. A {\it cubical object} in $\mA$ is an arbitrary
diagram $X: \Box^{op}\to \mA$.
In particular, {\it cubical sets} are diagrams $X: \Box^{op}\to \Set$,
and {\it cubical Abelian groups} are diagrams $F: \Box^{op}\to \Ab$.
A {\it cocubical object} is a diagram $\Box\to \mA$.
In particular, a {\it cocubical Abelian group} is a diagram of Abelian groups on $\Box$.

The category $\Box$ can be specified by the {\it conditional graph}.
This graph has vertices $\II^n$, where $n$ runs all nonnegative integers.
Its edges are the {face morphisms}
 $\II^{n-1}\stackrel{\delta^{n,\tau}_i}\to \II^{n}$
and degeneracy morphisms $\II^n\stackrel{\epsilon^n_i}\to \II^{n-1}$,
where
$n\geq 1$, $1\leq i\leq n$, $\tau\in \{0,1\}$.
Commutativity relations have been given above
(\ref{rel1})-(\ref{rel3}).
So, a cubical object $X: \Box^{op}\to \mA$ can be specified as a tuple  $(X_n, \partial^{n,\tau}_i, \sigma^n_i)$ consisting of sequences of objects $X_n=X(\II^n)$,
$n\geq 0$, and morphisms
$\partial^{n,\tau}_i=X(\delta^{n,\tau}_i): X_n\to X_{n-1}$,
$\sigma^n_i=X(\epsilon^n_i): X_{n-1}\to X_n$ satisfying to relations dual to (\ref{rel1})-(\ref{rel3}).
Morphisms $\partial^{n,\tau}_i$ are called {\it face operators}, and $\sigma^n_i$ are called
{\it degeneracy operators}.

\subsection
{Nondegenerate cubes}
Let $(X_k, \partial^{k,\tau}_i, \sigma^k_i)$ be a cubical set. Elements $x\in X_k$
are called  $k$-dimensional {\it cubes}.
The cube $x\in X_k$ is {\it degenerate} if there are
$y\in X_{k-1}$ and $i\in\{1, \ldots, k\}$ such that $\sigma^k_i(y)=x$.
Otherwise it is called  {\it nondegenerate}.

Let $h_{\II^n}: \Box^{op}\to \Set$ be the representable contravariant functor,
it is defined by $h_{\II^n}(\II^k)= \Box(\II^k,\II^n)$ on objects, and it assigns to each
morphism $f: \II^m\to \II^k$ the natural transformation
$\Box(f,\II^n): \Box(\II^k,\II^n)\to \Box(\II^m,\II^n)$ that converts
the morphism $g\in \Box(\II^k,\II^n)$ into $gf\in \Box(\II^m,\II^n)$.
The functor $h_{\II^n}$ is the cubical set and it is called the
{\it standard $n$-dimensional cube}. It can be given as the triple
$(X_k, \partial^{k,\tau}_i, \sigma^k_i)$ consisting of the sequence  $X_k=\Box(\II^k,\II^n)$
with maps
$\partial^{k,\tau}_i(f)=f\delta^{k,\tau}_i$ and
$\sigma^k_i(g)= g\epsilon^k_i$.

The following lemma follows directly from the definition of degenerate cube.
\begin{lemma}\label{degen}
A cube $f\in h_{\II^n}(\II^k)$ is degenerate if and only if
there exist a morphism
$g: \II^{k-1}\to \II^n$ and a number $i$, $1\leq i\leq k$, satisfying $f= g\epsilon^k_i$.
\end{lemma}

\begin{propos}\label{nondeg}
Let $n$ and $k$ be nonnegative integers. A cube
 $f\in h_{\II^n}(\II^k)$ is nondegenerate if and only if a corresponding map
 $f: \II^k\to \II^n$ is injective.
\end{propos}
{\sc Proof:}
Formula (\ref{grandis}) shows that each morphism $f: \II^k\to \II^n$
admits a decomposition $f=\delta\epsilon$ into a composition of a surjection and
an injection
$\II^k\stackrel{\epsilon}\to \II^{k-r}
\stackrel{\delta}\to \II^n$, where
$\epsilon= \epsilon^{k-r+1}_{i_1}\cdots \epsilon^k_{i_r}$, for some $i_1<\cdots< i_r$.
If $i<j+1$, then the equality (\ref{rel2}) shows that
$\epsilon^{k-1}_i\epsilon^k_{j+1}=\epsilon^{k-1}_j\epsilon^k_i$. It follows that
for each $q$ from the range $1\leq q\leq r$, there is a permutation
carrying  $\epsilon^k_{i_q}$ to the last place. Consequently, each morphism
that is not an injection, admits a decomposition
$\II^k\stackrel{\epsilon^k_i}\to \II^{k-1}\stackrel{g}\to \II^n$. So, it is the degenerate cube
of the standard cube $h_{\II^n}$ by Lemma \ref{degen}.
Obviously the inverse statement: if such a decomposition exists, then the morphism is not
an injection.
We obtain that the set of degenerate cubes
$\Box(\II^{k-1},\II^n)\epsilon^k_1\cup \cdots \cup \Box(\II^{k-1},\II^n)\epsilon^k_k$
contains all $f\in\Box(\II^k,\II^n)$ which are not injective.
\hfill$\Box$

\subsection{The construction of a projective resolution}

In this subsection, we construct a projective resolution of the cocubical
Abelian group $\Delta_{\Box}\ZZ$  in the category $\Ab^{\Box}$.

Consider the sequence of objects and morphisms in the category $\Ab^{\Box}$
\begin{equation}\label{complex}
0 \stackrel{d_0}\leftarrow \ZZ h^{\II^0} \stackrel{d_1}\leftarrow
\ZZ h^{\II^1} \stackrel{d_2}\leftarrow \ZZ h^{\II^2} \stackrel{d_3}\leftarrow
\ZZ h^{\II^3}
\leftarrow \cdots ,
\end{equation}
consisting of the cocubical groups and the natural transformations
$d_k= \sum\limits^k_{i=1}(-1)^i(\partial^{k,0}_i-\partial^{k,1}_i)$.
Here $\partial^{k,\tau}_i: \ZZ h^{\II^k}\to \ZZ h^{\II^{k-1}}$ are the natural transformations
whose components
$(\partial^{k,\tau}_i)_{\II^n}:\ZZ\Box(\II^k,\II^n) \to \ZZ\Box(\II^{k-1},\II^n)$
are defined on objects $\II^n\in \Box$ as taking values on the basis elements $f\in \Box(\II^k,\II^n)$
by means of equality $(\partial^{k,\tau}_i)_{\II^n}(f)= f\delta^{k,\tau}_i$.
From the relations $\partial^{k-1,\alpha}_i\partial^{k,\beta}_j=
\partial^{k-1,\beta}_{j-1}\partial^{k,\alpha_i}_j$ following from the equality (\ref{rel1})
we will follow the fulfillment of the formula $d_k d_{k+1}=0$, for all $k\geq 0$.
Hence, the sequence (\ref{complex}) is a chain complex.
But it is not exact in dimensions $k>0$ \cite[Remark 4.4]{X2008} and therefore can not be a resolution
of some object in $\Ab^{\Box}$.

For $n\geq 0$ and $k\geq 0$, consider a subset
$D_k(\II^n)\subseteq \Box(\II^k, \II^n)$ consisting of degenerate cubes
of the cubical set $h_{\II^n}$. The morphisms
$\Box(\II^k, \delta^{n,\tau}_i)$ and $\Box(\II^k, \epsilon^{n}_i)$
carry degenerate cubes into degenerate cubes. Consequently, $D_k$ is the subfunctor
of the functor $h^{\II^k}$.
Consider the embedding of the functors $\ZZ D_k\subseteq
\ZZ h^{\II^k}$, for $k\geq 0$. Its cokernel is the functor $\ZZ h^{\II^k}/\ZZ D_k$
defined on objects as the quotient groups $\ZZ h^{\II^k}(\II^n)/\ZZ D_k(\II^n)$. The
projection
 $\ZZ h^{\II^k}\to \ZZ h^{\II^k}/\ZZ D_k$ consists of the projections to the guotient groups.

 We have the following exact sequence in $\Ab^{\Box}$:
\begin{equation}\label{factor}
0 \leftarrow \ZZ h^{\II^k}/\ZZ D_k \stackrel{pr}\leftarrow \ZZ h^{\II^k}
\stackrel{\supseteq}\leftarrow \ZZ D_k \leftarrow 0.
\end{equation}
\begin{lemma}\label{proj}
The functor $\ZZ h^{\II^k}/\ZZ D_k$ is the projective object in the category
 $\Ab^{\Box}$.
\end{lemma}
{\sc Proof:} Let's construct a natural transformation $r: \ZZ h^{\II^k}\to \ZZ D_k$, 
the inverse from the left to the embedding $\ZZ D_k\subseteq\ZZ h^{\II^k}$.
By Yoneda's lemma, in order to construct $r$, it is enough
to specify the element $z= r_{\II^k}(1_{\II^k})\in \ZZ D_k(\II^k)$. 
And then the natural transformation $r$ will have components $r_{\II^n} (\alpha)= \alpha\circ z$,
 for all $n\geq 0$ and $\alpha\in\Box(\II^k, \II^n)$.
Resorting to the idea of Eilenberg and Mac\,Lane, used
in \cite[Prop. 7.2]{eil1953} when proving the representability
of normalized singular cubic homology groups, we put
$$
z = 1 -(1-\delta^{k,0}_1\sigma_1)(1- \delta^{k,0}_2\sigma_2)
\cdots(1- \delta^{k,0}_k\sigma_k).
$$
This $z$ is equal to the linear combination of products
$\delta^{k,0}_{s_1}\sigma_{s_1} \delta^{k,0}_{s_2} \sigma_{s_2}\cdots
\delta^{k,0}_{s_m}\sigma_{s_m}$ for some $1\leq m\leq k$ and $s_1 < \ldots <s_m$.

Each of these products is not a monomorphism, because multiplying this product on the right by $(1-\delta^{k,0}_{s_m}\sigma_{s_m})$ gives $0$. So $z$ is equal to a linear combination of degenerate morphisms, and $z\in\ZZ D_k(\II^k)$.
Hence, for any $n\geq 0$ and $\alpha\in\ZZ h^{\II^k}(\II^n)$ element $r_{\II^n}(\alpha)$ belongs 
to $\ZZ D_k(\II^n)$.

For any $1\leq i < j\leq k$, the operations $\delta^{k,0}_i\sigma^k_i$
and $\delta^{k,0}_j\sigma^k_j$ are permutable.
If $\alpha\in D^k(\II^n)$, then $\alpha= x\sigma_i$, for some
$1\leq i\leq k$ and $x\in h^{\II^k}(\II^n)$.
By virtue of the observed permutation property, 
$z= 1 - (1- \delta^{k,0}_i\sigma_i)y$ for some $y\in\Box(\II^k, \II^k)$.
From where
$$
x\sigma^k_i\circ z= x\sigma^k_i(1 - (1- \delta^{k,0}_i\sigma^k_i)y)=
x\sigma^k_i,
$$
and means $r_{\II^n}(\alpha)= \alpha$ for all $\alpha\in\ZZ D_k(\II^n)$.

Hence, $r$ is a retraction of $\ZZ h^{\II^k}$ by $\ZZ D_k$, and
the short exact sequence (\ref{factor}) splits. Since $\ZZ h^{\II^k}$ is a projective object in $\Ab^{\Box}$,
then $\ZZ h^{\II^k}/\ZZ D_k$ is projective.

The section $s: \ZZ h^{\II^k}/ \ZZ D_k\to \ZZ h^{\II^k}$ of the natural transformation $pr$ is determined by $r: \ZZ h^{\II^k}\to \ZZ D_k$ 
standard \cite[I.4]{mac1963}, according to the formula
$$
s_{\II^n}(\alpha+\ZZ D_k(\II^n))= \alpha -\alpha z.
$$
\hfill$\Box$

\smallskip
We proceed to the next step on constructing the projective resolution of the object
 $\Delta_{\Box}\ZZ\in \Ab^{\Box}$.
For this purpose we consider the homomorphisms
$(d_k)_{\II^n}: \ZZ h^{\II^k}(\II^n)\to \ZZ h^{\II^{k-1}}(\II^n)$
assigning to every element
$f:\II^k\to \II^n$ the sum
$\sum\limits_{i=1}^n (-1)^i(f\circ\delta^{k,0}_i-f\circ\delta^{k,1}_i)$.
For any $f\in D_k(\II^n)$, for $k\geq 1$,
there are a morphism $g: \II^{k-1}\to \II^n$
and a number $j\in \{1, \cdots, k\}$ such that $f=g\circ\epsilon^k_j$.
For $1\leq i\leq k$, we have the equality $f\delta^{k,\tau}_i=g\epsilon^k_j\delta^{k,\tau}_i$.
By (\ref{rel3}), for $i<j$, we have
$g\epsilon^k_j\delta^{k,\tau}_i=
g\delta^{k-1,\tau}_i\epsilon^{k-1}_{j-1}\in D_{k-1}(\II^n)$.
Similarly, for $i>j$, the equality $g\epsilon^k_j\delta^{k,\tau}_i\in D_{k-1}(\II^n)$ holds.
If $i=j$, then $g\epsilon^k_j\delta^{k,\tau}_i=g$, and the same time,
$f\circ\delta^{k,0}_i-f\circ\delta^{k,1}_i=g-g=0$. It follows that the homomorphisms
 $(d_k)_{\II^n}$ carry elements of
$\ZZ D_k(\II^n)$ into elements of $\ZZ D_{k-1}(\II^n)$.

So, we have the chain complex consisting of the cocubical Abelian groups $\ZZ h^{\II^k}/\ZZ D_k$
and the differentials
 $\overline{d}_k$ with components defined on the cosets of the subgroups
$\ZZ D_k(\II^n) \subseteq \ZZ h^{\II^k}(\II^n)$ by
$$
 (\overline{d}_k)_{\II^n}(f+\ZZ D_k(\II^n))= (d_k)_{\II^n}(f)+\ZZ D_{k-1}(\II^n).
$$

It follows from $D_0(\II^n)=\emptyset$ that $\ZZ h^{\II^0}/\ZZ D_0=\ZZ h^{\II^0}$.

\begin{lemma}\label{exactseq}
Let $K_*$ be the complex of Abelian groups
$$
 0 \leftarrow \ZZ h^{\II^0}/ \ZZ D_0(\II^n)
 \stackrel{(\overline{d}_1)_{\II^n}}\leftarrow
\ZZ h^{\II^1}/ \ZZ D_1(\II^n)
 \stackrel{(\overline{d}_2)_{\II^n}}\leftarrow
\ZZ h^{\II^2}/ \ZZ D_2(\II^n)
 \stackrel{(\overline{d}_3)_{\II^n}}\leftarrow \cdots.
$$
Then the homology groups  $H_q(K_*)$ equal $0$ for  $q>0$, and $H_0(K_*)=\ZZ$.
\end{lemma}
{\sc Proof:} We want to prove that the complex $K_*$ is isomorphic to the complex $C_*$
constructed in \cite[Lemma 4.1]{X2008} and consisting of the Abelian groups
and homomorphisms
$$
 0\leftarrow \ZZ\Box_+(\II^0,\II^n) \stackrel{d^+_1}\leftarrow
\ZZ\Box_+(\II^1,\II^n) \stackrel{d^+_2}\leftarrow \cdots
\stackrel{d^+_n}\leftarrow \ZZ\Box_+(\II^n,\II^n)\leftarrow 0.
$$
Here $d^+_k= \sum\limits^k_{i=1}(\ZZ(\partial^{k,1}_i)-\ZZ(\partial^{k,0}_i))$,
where $\partial^{k,\alpha}_i(\nu)=\Box_+(\delta^{k,\alpha}_i, \II^n)(\nu)= \nu\circ\delta^{k,\alpha}_i$
for all $\nu\in \Box_+(\II^k, \II^n)$ and $\alpha\in \{0,1\}$.

For this purpose we consider the homomorphisms $\gamma_k: C_k\to K_k$,
defined by $\gamma_k(\mu)= \mu + \ZZ D_k (\II^n)$,
where $\mu\in \ZZ\Box_+(\II^k,\II^n)$ is an arbitrary
linear combination  of injective morphisms $\II^k\to \II^n$ with integer
coefficients.
It is easy to see that the sequence of $\gamma_k$, $k\geq 0$, is an isomorphism of
the complexes $C_*\to K_*$.
According to \cite[Lemma 4.1]{X2008}, the homology groups of $C_*$
are equal $0$ in dimensions $k>0$ and $H_0(C_*)=\ZZ$.
Consequently, the same is true for $K_*$.
\hfill$\Box$

Deline the natural transformation
$\epsilon: \ZZ h^{\II^0}\to \Delta_{\Box}\ZZ$
as having the components  $\epsilon_{\II^n}: \ZZ\Box(\II^0,\II^n)\to \ZZ$
taking the values $\epsilon_{\II^n}(x)=1$ on basis elements $x\in \Box(\II^0, \II^n)$.

\begin{propos}\label{normres}
The sequence of the objects and morphisms
$$
0 \leftarrow \Delta_{\Box}\ZZ \stackrel{\epsilon}\leftarrow \ZZ h^{\II^0}/\ZZ D_0
\stackrel{\overline{d}_1}\leftarrow \ZZ h^{\II^1}/\ZZ D_1
\stackrel{\overline{d}_2}\leftarrow \ZZ h^{\II^2}/\ZZ D_2\leftarrow \cdots
$$
in the category
$\Ab^{\Box}$ is the projective resolution of the diagram $\Delta_{\Box}\ZZ$.
\end{propos}
{\sc Proof:} It follows from Lemma \ref{exactseq} that this sequence is exact.
 The cocubical Abelian groups $\ZZ h^{\II^k}/ \ZZ D_k$ are projective objects
 of the category $\Ab^{\Box}$ by Lemma \ref{proj}.
Therefore, this sequence is the projective resolution.
\hfill$\Box$

\subsection{Homology of the normalized complex of a cubical Abelian group}

Let  $F: \Box^{op}\to \Ab$ be a cubical Abelian group.
We construct its normalized complex $(C^N_k(F), d^N_k)$.

For any category $\mA$ with finite coproducts, its object
 $A\in \mA$ and a nonnegative integer $n\geq 0$, we denote by
 $A^{\sqcup n}$ the coproduct of $n$ copies $A\coprod \cdots \coprod A$
 of the object $A$. Let
 $in_i: A\to A^{\sqcup n}$ be the morphisms of the coproduct cocone,
 $1\leq i\leq n$.
 For any objects $A$, $B$ and morphisms $f_1, \ldots, f_n \in \mA(A,B)$
of the category $\mA$, we denote by
 $(f_1, \ldots, f_n): A^{\sqcup n}\to B$ the morphism such that
$(f_1, \ldots, f_n)\circ in_i= f_i$ for all $1\leq i\leq n$.

If  $\mA$ is an Abelian category, then the coproduct is denoted by the symbol $\oplus$.
In this case, for any
$f_1, \cdots, f_n \in \mA(A, B)$ there is the cokernel of the morphism
$(f_1, \ldots, f_n): A^{\oplus n}\to B$.
It consists of an object  $coker(f_1, \ldots, f_n)$ and a projection
  $pr: B\to coker(f_1, \ldots, f_n)$.

In the Abelian category, for a morphism  $f: A\to B$, the  object $coker(f)$
is defined up to isomorphism.
But if $\mA=\Ab$ then
as a cokernel, one can take a quotient group together with the
projection $pr: B\to B/ Im(f)$ assigning to every
$b\in B$ its coset $b+Im(f)$.
This cokernel we will call {\it canonical}.
The canonical cokernels give a functor on the category of
morphisms of Abelian groups.
Similarly, we can define a canonical cokernel in the category of diagrams  $\Ab^{\mD}$
on a small category  $\mD$.

The canonical cokernel of the homomorphism
$(f_1, \ldots, f_n): A^{\oplus n}\to B$ is equal to the quotient group $B/(Im(f_1) + \ldots + Im(f_n))$
with the projection assigning to every $b\in B$ its coset
$b+Im(f_1) + \ldots + Im(f_n)$.

Let $F\in \Ab^{\Box^{op}}$ be a cubical Abelian group.
Consider the chain complex  $(C_k(F),d_k)$ of the Abelian groups $C_k(F)=F(\II^k)$
and the differentials $d_k=\sum\limits^k_{i=1}(-1)^i (F(\delta^{k,0}_i)-F(\delta^{k,1}_i))$.
Let $D_k(F)= Im(F(\epsilon^k_1))+\cdots +Im(F(\epsilon^k_k)) \subseteq F(\II^k)$ be
the subgroup generates by the images of homomorphisms $F(\epsilon^k_i): F(\II^{k-1})\to F(\II^k)$,
for $1\leq i\leq k$. Its elements are called {\it degenerate chains}.

\begin{lemma}\label{degtodeg}
For each $k\geq 1$, the differential $d_k: C_k(F)\to C_{k-1}(F)$ carries degenerate chains into
degenerate chains.
\end{lemma}
{\sc Proof:} The differential $d_k$ carries each element $F(\epsilon^k_j)(x)\in Im(F(\epsilon^k_j))$
into a sum of the terms
$$
(-1)^i(F(\delta^{k,0}_i)-F(\delta^{k,1}_i))F(\epsilon^k_j)(x)= (-1)^i(F(\epsilon^k_j\delta^{k,0}_i)-F(\epsilon^k_j\delta^{k,1}_i))(x).
$$
For $i=j$ the term equals $0$. For $1\leq i<j\leq k$ it is equal to
 $$
 (-1)^i(F(\delta^{k-1,0}_i\epsilon^{k-1}_{j-1})-F(\delta^{k-1,1}_i\epsilon^{k-1}_{j-1}))(x)
\in Im F(\epsilon^{k-1}_{j-1}).
$$
 For $1\leq j<i \leq k$ it is equal to
$$
(-1)^i(F(\delta^{k-1,0}_{i-1}\epsilon^{k-1}_{j})-F(\delta^{k-1,1}_{i-1}\epsilon^{k-1}_{j}))(x)
\in Im F(\epsilon^{k-1}_{j}).
$$
Hence,  $d_k(D_k(F))\subseteq D_{k-1}(F)$.
\hfill$\Box$

\begin{definition}\label{cnk}
The {\it normalized complex} of a cubical Abelian group $F\in \Ab^{\Box^{op}}$
consists of the quotient group $C^N_k(F)= F({\II}^k)/D_k(F)$
and differentials $d^N_k: C^N_k(F)\to C^N_{k-1}(F)$ defined by $d^N_k(f+D_k(F))= d_k(f)+ D_{k-1}(F)$ for all integers $k\geq 0$. It is defined by
$C^N_k(F)=0$ if $k<0$.
\end{definition}

Let $\Comp$ be the category of chain complexes of Abelian groups $\{(K_n, d_n)\}$ such that
$K_n=0$ for all $n<0$.
It follows from Lemma \ref{degtodeg} that we can define the differentials
$d'_k: D_k(F)\to D_{k-1}(F)$ by $d'_k(f)= d_k(f)$ for all $k\geq 1$.
\begin{propos}\label{functorchain}
The maps assigning to $F\in \Ab^{\Box^{op}}$
the chain complexes  $C_*(F)=(C_k(F), d_k)$, $D_*(F)=(D_k(F), d'_k)$,
and $C^N_*(F)=(C^N_k(F), d^N_k)$ can be considered as functors $\Ab^{\Box^{op}}\to \Comp$.
\end{propos}
{\sc Proof:} The functor $C_*$ is defined in the usual way.
It follows by Lemma \ref{degtodeg} that $D_*$ is a functor.
By Definition \ref{cnk}, we conclude that $C^N_*$ can be considered as a functor.
\hfill$\Box$

Image of the homomorphism $F(\II^{k-1})^{\oplus k}
\xrightarrow
{(F(\epsilon^k_1,\cdots, F(\epsilon^k_k)))}
F(\II^k)$ equals $D_k(F)$, hence we have from Definition
 \ref{cnk} the following exact sequence:
\begin{equation}\label{cnkrepr}
F(\II^{k-1})^{\oplus k}
\xrightarrow
{(F(\epsilon^k_1,\cdots, F(\epsilon^k_k)))}
F(\II^k)\to C^N_k(F)\to 0.
\end{equation}

Denote $P_k = \ZZ h^{\II^k}/ \ZZ D_k$.
Let $P_*=(P_k, \overline{d}_k)$ be the projective resolution constructing in
Proposition \ref{normres}.

\begin{theorem}\label{main1}
The complex  $P_*\otimes F$ is isomorphic to the complex $C^N_*(F)$.
\end{theorem}
{\sc Proof:} The cubical Abelian groups  $F$ and $\ZZ h^{\II^{(-)}}\otimes F$ are isomorphic.
The natural isomorphism is given on the Fig. 1.
The application of Proposition \ref{functorchain} leads to the isomorphism of
complexes
$C^N_*(F)\cong C^N_*(\ZZ h^{\II^{(-)}}\otimes F)$.

Now it is suficient to prove an isomorphism of complexes $C^N_*(\ZZ h^{\II^{(-)}}\otimes F)$
and $P_*\otimes F$.
First, we constructing isomorphisms  $C^N_k(\ZZ h^{\II^{(-)}}\otimes F)\to P_k\otimes F$,
for all $k\geq 0$.
By Definition \ref{cnk}, we have
$C^N_k(\ZZ h^{\II^{(-)}}\otimes F)= (\ZZ h^{\II^k}\otimes F) / D_k (\ZZ h^{\II^{(-)}}\otimes F)$.
Hence, we need an isomorphism
$(\ZZ h^{\II^k}\otimes F) / D_k (\ZZ h^{\II^{(-)}}\otimes F)\to (\ZZ h^{\II^k}/\ZZ D_k)\otimes F$.

Consider the exact sequence constructed in the same way as
(\ref{cnkrepr}):

\begin{equation}\label{rightexact}
(\ZZ h^{\II^{k-1}})^{\oplus k}
\xrightarrow{(\ZZ h^{\epsilon^k_1}, \cdots, \ZZ h^{\epsilon^k_k})} \ZZ h^{\II^k}\to P_k \to 0.
\end{equation}
The functor $(-)\otimes F: \Ab^{\Box}\to \Ab$ is right exact, so it
carries the exact sequence (\ref{rightexact}) into the exact sequence
$$
(\ZZ h^{\II^{k-1}})^{\oplus k}\otimes F
\xrightarrow{(\ZZ h^{\epsilon^k_1}, \cdots,
\ZZ h^{\epsilon^k_k})\otimes 1_F}
\ZZ h^{\II^k}\otimes F \to
P_k \otimes F \to 0.
$$
The first homomorphism of this sequence admits the following decomposition
$$
(\ZZ h^{\II^{k-1}})^{\oplus k}\otimes F
\stackrel{e_k\otimes F}\to
\ZZ D_k \otimes F
\stackrel{j_k\otimes 1_F}\to \ZZ h^{\II^k} \otimes F ,
$$
where $j_k$ is denoted the inclusion $\ZZ D_k\subseteq \ZZ h^{\II^k}$, and $e_k$ is
the epimorphism onto $\ZZ D_k$.
There is an exact sequence of Abelian groups
$Tor(P_k, F)\to \ZZ D_k \otimes F
\stackrel{j_k\otimes 1_F}\to \ZZ h^{\II^k} \otimes F \to P_k\otimes F \to 0$.
By Lemma \ref{proj}, the object $P_k$ is projective.
Therefore, $j_k\otimes F$ is the monomorphism, and $e_k\otimes 1_F$
with $j_k\otimes 1_F$ form the decomposition of  $(\ZZ h^{\epsilon^k_1}, \cdots,
\ZZ h^{\epsilon^k_k})\otimes 1_F$  into a composition of an
epimorphism and a monomorphism.

The finite direct sums commute with the tensor product. Hence, there is
an isomorphism
$(\ZZ h^{\II^{k-1}})^{\oplus k}\otimes F
\stackrel{\alpha_k}\to (\ZZ h^{\II^{k-1}}\otimes F)^{\oplus k}$.
It leads to the commutative diagram shown in Fig. 2.
\begin{figure}[h]
$$
\xymatrix{
(\ZZ h^{\II^{k-1}})^{\oplus k}\otimes F\ar@{->>}[d]_{e_k\otimes 1_F}
\ar[rr]^{\alpha_k}_{\cong}
 && (\ZZ h^{\II^{k-1}}\otimes F)^{\oplus k}\ar@{->>}[d]\\
\ZZ D_k \otimes F \ar@{-->}[rr]_{\cong}^{\overline{\alpha}_k}
\ar@{>->}[rd]_{j_k\otimes 1_F} && D_k(\ZZ h^{\II^{(-)}}\otimes F)\ar[ld]^{\subseteq}\\
 & \ZZ h^{\II^k}\otimes F
}
$$
\centerline{\small Fig. 2. The isomorphism
$\ZZ D_k \otimes F\to D_k(\ZZ h^{\II^{(-)}}\otimes F)$.}
\end{figure}

Any two decompositions by an epimorphism and a monomorphism are isomorphic.
It follows that there exists a unique isomorphism $\overline{\alpha}_k$, shown by the pointed arrow in Fig.2 making a commutative triangle and the square of the diagram.

Consider the isomorphism of the exact sequences shown in Fig. 3.
\begin{figure}[h]
$$
\xymatrix{
0 \ar[r] &\ZZ D_k\otimes F\ar[r]^{j_k\otimes 1_F} \ar[d]_{\overline{\alpha}_k}
    & \ZZ h^{\II^k}\otimes F\ar[r] \ar[d]^{=} & P_k\otimes F \ar[r] \ar[d]^{\cong}_{\beta_k}& 0\\
0 \ar[r] & D_k(\ZZ h^{\II^{(-)}}\otimes F) \ar[r]^(.6){i_k}_(.6){\subseteq}
    & \ZZ h^{\II^k}\otimes F\ar[r] & C^N_k(\ZZ h^{\II^{(-)}}\otimes F) \ar[r]& 0\\
}
$$
\centerline{\small Fig.3. The isomorphism of the exact sequences.}
\end{figure}
Here the isomorphism $\beta_k$ appears as an isomorphism of cokernels complementing the
diagram to a commutative diagram.
It remains for us to prove that the isomorphisms $\beta_k$ commute with the differentials.
Using the cokernel construction, we can obtain the remaining morphisms of the diagram
shown on Fig. 4.
\begin{figure}[h]
$$
\xymatrix{
P_k\otimes F\ar[d]_{d_k\otimes 1_F} \ar[r]^(.4){\cong}_(.4){\beta_k}
& C^N_k(\ZZ h^{\II^{(-)}}\otimes F)
\ar[d]^{d^N_k}\\
P_{k-1}\otimes F \ar[r]_(.4){\cong}^(.4){\beta_{k-1}}
& C^N_{k-1}(\ZZ h^{\II^{(-)}}\otimes F)\\
}
$$
\centerline{\small Fig. 4. The isomorphism of the complexes
$P_*\otimes F\to C^N_*(\ZZ h^{\II^{(-)}}\otimes F)$.}
\end{figure}
This diagram is commutative as consisting of the cokernels of the morphisms
 $j_k\otimes 1_F$,  $j_{k-1}\otimes 1_F$, $i_k$ and $i_{k-1}$ of the commutative diagram
 shown on Fig. 5.
\begin{figure}[h]
$$
\xymatrix{
\ZZ D_k\otimes F
\ar[ddd]^{d'_k\otimes 1_F}
\ar[rd]^{j_k\otimes F}
\ar[rrr]^{\overline{\alpha}_k}
&&& D_k(\ZZ h^{\II^{(-)}}\otimes F)
\ar[ld]^{i_k} \ar[ddd]^{d_k}\\
& \ZZ h^{\II^k}\otimes F \ar[d]^{d_k\otimes 1_F}
\ar[r]^{=}
& \ZZ h^{\II^k}\otimes F\ar[d]^{d_k\otimes 1_F}\\
& \ZZ h^{\II^{k-1}}\otimes F \ar[r]^{=} & \ZZ h^{\II^{k-1}}\otimes F
\\
\ZZ D_{k-1}\otimes F
\ar[ru]^{j_{k-1}\otimes F}
\ar[rrr]^{\overline{\alpha}_{k-1}} &&&
D_{k-1}(\ZZ h^{\II^{(-)}}\otimes F)\ar[lu]^{i_{k-1}}
}
$$
\centerline{\small Fig. 5. To the proof of the isomorphism of complexes
$P_*\otimes F\to C^N_*(\ZZ h^{\II^{(-)}}\otimes F)$.}
\end{figure}
We arrive at an isomorphism of the complexes $P_*\otimes F\cong C^N_*(\ZZ h^{(-)}\otimes F)$.
The diagram on Fig. 1 gives an isomorphism $\ZZ h^{(-)}\otimes F \cong F$ in the category
$\Ab^{\Box^{op}}$ leading to  $C^N_*(\ZZ h^{(-)}\otimes F) \cong C^N_*(F)$. Therefore,
$P_*\otimes F\cong C^N_*(F)$.
\hfill$\Box$

\begin{corollary}\label{cor1}
For any cubical Abelian group $F$ and integer $n\geq 0$
the homology groups of the complex
$H_n(C^N_*(F))$ are isomorphic to groups $\coLim^{\Box^{op}}_n F$.
\end{corollary}
{\sc Proof:} It follows from Propositions \ref{projcolim}, \ref{normres} and Theorem \ref{main1}.
\hfill$\Box$

\begin{corollary}\label{split}
For any $k\geq 0$ and $F\in \Ab^{\Box^{op}}$  there is an isomorphism
$F(\II^k)\cong C^N_k(F)\oplus D_k(F)$.
\end{corollary}
{\sc Proof:} The normalized complex $C^N_*(F)$ is defined by the exact sequence
$0\to D_k(F)\to F(\II^k)\to C^N_k(F)\to 0$. The isomorphism of cubical Abelian groups
$F\cong \ZZ h^{\II^{(-)}}\otimes F$  leads to  the fact that this exact sequence is isomorphic to
$0\to D_k(\ZZ h^{\II^{(-)}}\otimes F) \to \ZZ h^{\II^{k}}\otimes F
\to C^N_k(\ZZ h^{\II^{(-)}}\otimes F)\to 0$. The last is isomorpic (Fig. 3)
to
$0\to \ZZ D_k\otimes F \stackrel{j_k\otimes 1_F}\to \ZZ h^{\II^k}\otimes F
\stackrel{pr_k\otimes 1_F}\to P_k\otimes F \to 0$.
The exact sequence
$0\to \ZZ D_k \stackrel{j_k}\to \ZZ h^{\II^k} \stackrel{pr_k}\to P_k\to 0$
is split. So, we can conclude that  the exact sequnce
$0\to \ZZ D_k\otimes F \stackrel{j_k\otimes 1_F}\to \ZZ h^{\II^k}\otimes F
\stackrel{pr_k\otimes 1_F}\to P_k\otimes F \to 0$ is split.
\hfill$\Box$

An assertion analogous to Corollary \ref{split} was obtained in
 \cite[Theorem 2.2]{leb2017} for skew cubical structures.

\section{Homology of cubical sets with coefficients in contravariants systems}

Let $X\in \Set^{\Box^{op}}$ be a cubical set. For integer $m\geq 0$, its $m$-dimensional
{\it singular cube} is an arbitrary morphism  $h_{\II^m}\stackrel{\xi}\to X$
of cubical sets.
Denote by $\Box/X$ the comma-category of objects $h^{\Box}$-over $X$
in the sense of \cite{mac2004} where  $h^{\Box}: \Box\to \Set^{\Box^{op}}$ is the Yoneda embedding. We follow to \cite{gab1967} and call it the {\it left fibre of $h^{\Box}$ over $X$}.
The objects of $\Box/X$ are singular cubes of  $X$.
By Yoneda Lemma, each singular cube $h_{\II^m}\stackrel{\xi}\to X$ is defined by the
element
$x=\xi_{\II^m}(1_{\II^m})\in X_m$, and we can denote it by $\xi=\widetilde{x}$.
So, these objects are equal to $h_{\II^m}\stackrel{\widetilde{x}}\to X$ for some
$x\in X_m$.
Morphisms of the category $\Box/X$ are commutative triangles shown in Fig. 6.
\begin{figure}[h]
$$
\xymatrix{
 h_{\II^m}\ar[rr]^{h_{\alpha}}\ar[rd]_{\widetilde{x}} && h_{\II^n}\ar[ld]^{\widetilde{y}}\\
 & X
}
$$
\centerline{\small Fig. 6. A morphism of $\Box/X$.}
\end{figure}
A {\it contravariant system} on $X$ is an arbitrary diagram of Abelian group
 $F: (\Box/X)^{op}\to \Ab$.
We define homology groups of a cubical set $X$ with coefficients in $F$ as the homology
groups of a normalized chain complex and we
prove that these homology groups are isomorphic to Abelian groups $\coLim^{(\Box/X)^{op}}_n F$
for all $n\geq 0$.

\subsection{Normalized complex of a cubical system with coefficients in a contravariant system}

Let $X$ be a cubical set.
The category $(\Box/X)^{op}$ can be considered as consisting
of the set of objects $x\in \coprod\limits_{n\geq 0}X_n$.
Its morphisms can be given as triples $x\stackrel{\alpha}\to y$ such that
$X(\alpha)(y)=x$.
The {\it forgetful functor} $Q_X: \Box/X \to \Box$ for the left fibre, assigns to every object $h_{\II^m}\stackrel{\widetilde{x}}\to X$ the object $\II^m$
and to every commutative triangle
 (Fig. 6) the morphism $\II^m\stackrel{\alpha}\to \II^n$.

For any diagram of Abelian groups $F: (\Box/X)^{op}\to \Ab$,
we denote by $(C_n(X,F), d^{n,\tau}_i, s^n_i)$ the cubical Abelian group consisting
of Abelian groups $C_n(X,F)=\bigoplus\limits_{x\in X_n}F(x)$ with face operators
$d^{n,\tau}_i:  C_{n}(X,F)\to C_{n-1}(X,F)$ and degeneracy operators
$s^n_i: C_{n-1}(X,F)\to C_{n}(X,F)$.
The operators $d^{n,\tau}_i$ are defined
as the homomorphisms that make commutative diagrams
shown on Fig. 7.
\begin{figure}[h]
$$
\xymatrix{
\bigoplus\limits_{x\in X_n}F(x) \ar[rrr]^{d^{n,\tau}_i} &&&
	\bigoplus\limits_{x\in X_{n-1}}F(x)\\
F(x)\ar[u]^{in_x}
\ar[rrr]_{F(\delta_i^{n,\tau}:x\to X(\delta^{n,\tau}_i)x)} &&&
F(X(\delta^{n,\tau}_i)x)\ar[u]_{in_{X(\delta^{n,\tau}_i)x}}
}
$$
\centerline{\small Fig. 7. The definition of face operators.}
\end{figure}
The operators  $s^n_i$ are defined
as the homomorphisms that make commutative diagrams
shown on Fig. 8.
\begin{figure}[h]
$$
\xymatrix{
\bigoplus\limits_{x\in X_{n-1}}F(x) \ar[rrr]^{s_i^{n}} &&&
	\bigoplus\limits_{x\in X_{n}}F(x)\\
F(x)\ar[u]^{in_x}
\ar[rrr]_{F(\epsilon_i^{n}:x\to X(\epsilon^{n}_i)x)} &&&
F(X(\epsilon^{n}_i)x)\ar[u]_{in_{X(\epsilon^n_i)x}}
}
$$
\centerline{\small Fig. 8. The definition of degeneracy operators.}
\end{figure}

Denote by $C^N_*(X,F)= (C^N_n(X,F),d^N_n)$ the
normalized complex of the cubical Abelian group
$(C_n(X,F), d^{n,\tau}_i, s^n_i)$. We call it
a {\it normalized complex of $X$ with coefficients in $F$}.
By Definition \ref{cnk}, it consists of quotient groups
 $C^N_n(X,F)= C_n(X,F)/D_n(X,F)$ where
$D_n(X,F)$ is a subgroup generated by images of the homomorphisms $s^n_i$.
Its differentials are defined by $d^N_n(a+D_n(X,F))= d_n(a)+D_{n-1}(X,F)$.

For an arbitrary functor between small categories $S: \mC\to \mD$
and a cocomplete category $\mA$, the functor $(-)\circ S: \mA^{\mD}\to \mA^{\mC}$
has a left adjoint functor  $\Lan^S: \mA^{\mC} \to \mA^{\mD}$ called
a {\it left Kan extension} \cite{mac2004}.
We apply this construction to the functor $S= Q_X^{op}: (\Box/X)^{op}\to \Box^{op}$
and $\mA= \Ab$. For any contravariant system  $F$ on $X$, it gives
the cubical Abelian group $\Lan^{Q_X^{op}}F$. It easy to see that $\Lan^{Q_X^{op}}F$
is isomorphic to the cubical Abelian group $(C_n(X,F), d^{n,\tau}_i, s^n_i)$.
See the proof in
\cite[Proposition 3.7]{X2008} for the general case of a small category
$\mD$ (instead $\Box$) and for a functor $F: (\mD/X)^{op}\to \Ab$.
It follows from Proposition \ref{functorchain} that the normalized complexes
$C^N_*(X,F)$ and $C^N_*(\Lan^{Q_X^{op}}F)$ are isomorphic.
\begin{definition}
Homology groups $H_n(X,F)$ of a cubical set $X$ with coefficients
in a contravariant system $F$ on $X$ are the homology groups
$H_n(C^N_*(X,F))$ of the normalized complex.
\end{definition}

\subsection{Main Theorem}

Before passing to the main theorem, we recall some definitions
and prove an auxiliary lemma.
For an arbitrary functor $S: \mC\to \mD$ and an object $d\in \mD$ the
(left) {\it fibre} $S/d$ of $S$ over $d$ \cite[Application 2]{gab1967} is
the category whose objects are pairs $(a\in \mC, \alpha: S(a)\to d)$
consisting of an object $a\in \mC$ and a morphism $\alpha: S(a)\to d$.
Its morphisms $(a_1,\alpha_2)\to (a_1, \alpha_2)$ are specified as morphisms $\gamma: a_1\to a_2$
such that $\alpha_2\circ S(\gamma)= \alpha_1$.
The {\it forgetful functor}  $Q_d: S/d\to \mC$ assigns to each pair $(a,\alpha)$
the object $a\in \mC$ and to any morphism $(a_1,\alpha_1)\stackrel{\gamma}\to (a_2,\alpha_2)$
the morphism $\gamma$.

\begin{lemma}\label{misc}
Let $\mD$ be a small category and $X \in \Set^{\mD^{op}}$ be a diagram of sets.
Then for any functor $F: (\mD/X)^{op}\to \Ab$
there is a natural isomorphism
$\coLim_n^{\mD^{op}}\Lan^{Q_X^{op}}F \cong \coLim_n^{(\mD/X)^{op}}F$.
\end{lemma}
{\sc Proof:} By \cite[Proposition 2, Remark 3.8]{gab1967}, the values of the left derived functors of
$\Lan^{Q_X^{op}}$ are equal to
$(\Lan_q^{Q_X^{op}}F)(d)= \coLim_q^{Q_X^{op}/d}F\circ Q_d^{op}$,
for all $d\in \mD$ and $q\geq 0$. Here $Q_d^{op}: Q_X^{op}/d \to (\mD/X)^{op}$
is the forgetful functor for fibre of $Q_X^{op}$ over $d$.
Every connected component of the category
$Q_X^{op}/d$ has a terminal object \cite[Example 1.1]{X1997} therefore
$\Lan_q^{Q_X^{op}}F= 0$, for all $q>0$. It follows that the spectral sequence of
 \cite[Application 2, Theorem 3.8]{gab1967} is degenerated into the isomorphisms
$ \coLim_n^{\mD^{op}}\Lan^{Q_X^{op}}F \cong \coLim_n^{(\mD/X)^{op}}F$.
\hfill$\Box$

\begin{theorem}\label{main2}
For each contravariant system $F$ on a cubical set $X$,
there are isomorphisms
$H_n(X,F)\cong \coLim^{(\Box/X)^{op}}_n F$ for all $n\geq 0$.
\end{theorem}
{\sc Proof:} The complex $C^N_*(X,F)$ is isomorphic to  $C^N_*(\Lan^{Q_X^{op}}F)$,
hence $H_n(X,F)\cong \coLim_n^{\Box^{op}}\Lan^{Q_X^{op}}F$ by Corollary \ref{cor1}.
The application of Lemma \ref{misc} completes the proof.
\hfill$\Box$

We give one of the simplest applications of Theorem \ref{main2}.
A contravariant system $L: (\Box/X)^{op}\to \Ab$ on a cubical set  $X$ is called
a {\it local system} if it consists of isomorphisms.
It follows from \cite[Application 2, Proposition 4.4]{gab1967} for all $n\geq 0$,
that there are isomorphisms
$\coLim^{(\Box/X)^{op}}_n L \stackrel{\cong}\to \coLim^{\Box/X}_n L^{-1}$,
where $L^{-1}$ is constructed by inversion of homomorphisms of the local system $L$.
In particular, for the standard $k$-dimensional cube $X=h_{\II^k}$,
the category $\Box/h_{\II^k}$ has the terminal object
$\widetilde{1_{\II^k}}: h_{\II^k}\to h_{\II^k}$.
Hence, the colimit of  $L^{-1}$ is equal to $L^{-1}(\widetilde{1_{\II^k}})$, and the
functors of $\coLim^{\Box/h_{\II^k}}_n: \Ab^{\Box/h_{\II^k}}\to \Ab$ are equal to $0$ for $n>0$.
Applying Theorem \ref{main2}, immediately obtain the following statement:

\begin{corollary}
For any local system $L$ on the standard cube $h_{\II^k}$, we have
$$
H_n(h_{\II^k}, L)=
\begin{cases}
L(\widetilde{1_{\II^k}}), & \mbox{ if } n=0;\\
0, & \mbox{ if } n>0.
\end{cases}
$$
\end{corollary}

\section{Spectral sequences for the homology of cubical sets}

To compute the homology groups of simplicial, semisimplicial, and precubical sets, it
suffices to use a chain complex in which the differentials are linear combinations of
face operators. In the case of a cubical set, this complex must be normalized,
and it is not clear how to study the properties of homology groups of
cubical sets with coefficients in contravariant systems consisting of homomorphisms that may
not be isomorphisms.

In this section it is shown that in many cases Theorem \ref{main2}
allows us to solve this problem.
The results analogous to those obtained in \cite{X2008} for precubical sets are proved.

\subsection{The isomorphism criterion for the homology groups of cubical sets}

Consider a morphism of cubical sets $f: X\to Y$.
The functor $\Box/f: \Box/X\to \Box/Y$ assigns to any singular cube
$\widetilde{x}: h_{\II^n}\to X$ the singular cube
$f\circ\widetilde{x}: h_{\II^n}\to Y$.
Let $f^*: \Ab^{(\Box/Y)^{op}} \to \Ab^{(\Box/X)^{op}}$ be the fuctor
assigning to each functor
 $F: (\Box/Y)^{op}\to \Ab$ the functor
$F\circ(\Box/f)^{op}$.
There are canonical homomorphisms
$\coLim^{(\Box/X)^{op}}_n f^* (F)\to
\coLim^{(\Box/Y)^{op}}_n F$. The application of Theorem \ref{main2} leads
to homomorphisms $H_n(X,f^*(F))\to H_n(Y,F)$.

For a cube $y\in Y_n$ of a cubical set  $Y$, the
{\it inverse image} $\overleftarrow{f}({y})$ of the singular cube
 $h_{\II^n}\stackrel{\widetilde{y}}\to Y$ is the limit of the diagram
$X\stackrel{f}\to Y \stackrel{\widetilde{y}}\leftarrow h_{\II^n}$.
Denote by $f_y: \overleftarrow{f}({y})\to X$ the cone morphism of this limit to $X$.
It easy to see that there is an isomorphism of categories
$(\Box/f)/\widetilde{y})\cong \Box/\overleftarrow{f}(y)$.
Applying \cite[Theorem 2.3]{obe1968}, \cite[Lemma 3.1]{X2008}, and Theorem \ref{main2},
we arrive to the following statement:

\begin{corollary}
Let $f: X\to Y$ be a morphism of cubical sets.
Then the following properties of the morphism $f$ are equivalent:
\begin{enumerate}
\item For each $y\in \coprod_{k\geq 0} Y_k$
the groups $H_n(\overleftarrow{f}(y),\Delta\ZZ)$ equal to $0$ for all $n>0$,
and $H_0(\overleftarrow{f}(y),\Delta\ZZ)\cong \ZZ$.
\item The canonical homomorphisms of Abelian groups
$H_n(X,f^{*}F)\to H_n(Y,F)$ are isomorphisms for every functor
$F: (\Box/Y)^{op}\to Ab$.
\end{enumerate}
\end{corollary}

We give an example showing that even morphisms between standard cubes do not
preserve homology groups with coefficients in contravariant systems.

\begin{example}
The $0$-dimensional standard cube $h_{\II^0}$ is terminal object in $\Set^{\Box^{op}}$.
Consider the (unique) morphism
$f=h_{\epsilon^1_1}: h_{\II^1}\to h_{\II^0}$.
For the singular cube  $\widetilde{y}=h_{\epsilon^1_1}: h_{\II^1}\to h_{\II^0}$
its inverse image  ${\overleftarrow{f}}(y)$ is isomorphic to the product
$h_{\II^1}\times h_{\II^1}$.
It was shown in \cite[Remark 3.5]{jar2006} that the geometric realization of
$h_{\II^1}\times h_{\II^1}$ has the homotopy type of the wedge
of spheres $S^2\vee S^1$.
But we will not use this for calculating the homology groups of the cubical set
 $h_{\II^1}\times h_{\II^1}$.
The set of $n$-dimensional singular cubes of this product
consists of pairs of functions
$(\alpha,\beta)\in \Box(\II^n,\II^1)\times \Box(\II^n, \II^1)$.
For $n\geq 3$, it does not have nondegenerate cubes.
It has two nondegenerate $2$-dimensional cubes:
$(\alpha(x_1,x_2)=x_1, \beta(x_1,x_2)=x_2)$ and
$(\alpha(x_1,x_2)=x_2, \beta(x_1,x_2)=x_1)$.
In standard notations, they are equal to $(\epsilon^2_2, \epsilon^2_1)$ and
$(\epsilon^2_1, \epsilon^2_2)$.
It has five nondegenerate $1$-dimensional cubes
$(id,0)$, $(id,1)$, $(id,id)$, $(0,id)$, $(1, id)$, where $id: \II^1\to \II^1$ is the
identity map and
 $0: \II^1\to \II^1$ (resp. $1: \II^1\to \II^1$) are the maps taking constant values equal
 to
 $0$ (resp. $1$).
Four pairs $(0,0)$, $(1,0)$, $(1,1)$, $(0,1)$ of the maps $\II^0\to \II^1$ are
nondegenerate as $0$-dimensional cubes.
Writing down the matrices of differentials and reducing them to the normal form
of Smith, we obtain
$H_0(h_{\II^1}\times h_{\II^1})\cong H_1(h_{\II^1}\times h_{\II^1})\cong
 H_2(h_{\II^1}\times h_{\II^1}) \cong\ZZ$.
Here $H_n(X)$ are the homology groups of cubical sets with
coefficients in $\Delta\ZZ$.
In particular,
 there exists a contravariant system $F$ on the point $h_{\II^0}$,
such that for $f=h_{\epsilon^1_1}: h_{\II^1}\to h_{\II^0}$ the morphism $H_n(h_{\II^1},f^* F)\to H_n(h_{\II^0},F)$ is not an isomorphism.
\end{example}

\subsection{Homology groups of the colimits of cubical sets}

Let  $X: J\to \mA$ be a diagram of objects in a category $\mA$
defined on a small category $J$. In some cases, it is convenient to denote it
by
 $\{X^i\}_{i\in J}$ (shortly $\{X^i\}$),
 indicating its values
 $X^i$ on objects $i\in \Ob(J)$.
 Denote by $\lambda_i: X^i\rightarrow \coLim^J\{X^i\}_{i\in J}$ the morphisms of the colimit cone.

We will consider the first quadrant spectral sequences in the sense of \cite{mac1963}.
The following statement is proved just like
\cite[Proposition 5.2]{X2008}.
But we need to use Theorem \ref{main2} instead of \cite[Theorem 4.3]{X2008}.

\begin{corollary}\label{covering}
Let $J$ be a small category and let $\{X^i\}_{i\in J}$ be a diagram of cubical sets such that
\begin{equation}\label{c3}
\coLim_q^J\{\ZZ(X^i_n)\}_{i\in J}= 0, \mbox{ for all } n\in \NN \mbox{ and }
q>0.
\end{equation}
Then for any contravariant system
$F$ on $\coLim^J\{X^i\}$
there exists a first quadrant spectral sequence
$$E^2_{p,q}= \coLim_p^{J}\{H_q(X^i,\lambda_i^*(F))\}_{i\in J} \Rightarrow
\{H_{p+q}(\coLim^J\{X^i\}_{i\in J}, F)\}.
$$
\end{corollary}

Every monoid $J$ can be considered as the category with a unique
object and the set of morphisms $J$.
 Any left $J$-set $X$ defines the functor ${\cal X}: J\to \Set $
 assigning to its unique object the set $X$ and to every morphism
  $g\in J$ the map defined by ${\cal X}(g)(x)= g x$.

The monoid $J$ {\it acts free} on $X$,
 if the functor ${\cal X}$ is isomorphic to the coproduct of representable functors
  $J\to \Set$.

The monoid $J$ {\it acts free on a cubical set} $X$
if it acts free on $X_n$ for each $n\geq 0$.

If $J$ is a monoid, then the values of left
derived functors of the colimit $\coLim^J_n A$ are equal to the homology groups
${\cal H}_n(J, A)$ of $J$ with coefficients in $J$-module $A$.
If a monoid $J$ acts free on a cubical set $X$, then the condition (\ref{c3}) is satisfied.

\begin{corollary}
Let a monoid $J$ acts free on a cubical set $X$, and let
$\lambda: X\to X_J$ be the canonical projection onto
 quotient of the cubical set $X$ by $J$-action.
 Then for every contravariant system
$F$ on $X_J$ there is a spectral sequence
$
E^2_{p,q}= {\cal H}_p(J, H_q(X,\lambda^*(F))) \Rightarrow
H_{p+q}(X_J, F).
$
\end{corollary}

A {\it locally directed covering of a cubical set $X$} is a diagram of cubical sets $\{X^i\}_{i\in J}$
on a partially ordered set
$J$ satisfying the following conditions:

1) this diagram assign to every pair of elements $i\leq j$ of $J$
the inclusion $X^i \subseteq X^j$;

2) for each $n\in \NN$ there is equality $X_n= \bigcup_{i\in J}X_n^i$;

3) for every  $\sigma\in X_n^i\cap X_n^j$ there is an element $k\in J$ such that
$k\leq i$, $k\leq j$ and $\sigma\in X_n^k$.

In this case the conditions of Corollary \ref{covering} are satisfied and we obtain
the following statement:

\begin{corollary}\label{locdir}
Let $\{X^i\}_{i\in J}$ be a locally directed covering of a cubical set
$X$ and let
$\lambda_i$ be inclusions $X^i \subseteq X$.
Then for any contravariant system
$F$ on $X$, there exists a first quadrant spectral sequence
$
E^2_{p,q}= \coLim_p^{J}\{H_q(X^i,\lambda_i^*(F))\}_{i\in J}
 \Rightarrow  \{H_{p+q}(X, F)\}.
$
\end{corollary}

In particular, this implies the existence of an exact sequence for the union of cubical sets:
\begin{corollary}
For any contravariant system $F$ on the union of cubical sets $X^1\cup X^2$,
there is the exact sequence
\begin{multline*}
\cdots \to H_{n+1}(X^1\cup X^2, F) \to H_n(X^1\cap X^2, \lambda^*_0 F)\\
 \to
H_n(X^1, \lambda^*_1 F)\oplus H_n(X^2, \lambda^*_2 F) \to
H_n(X^1\cup X^2, F) \to \cdots
\end{multline*}
for all $n\geq 0$.
Here $\lambda_i: X^i\to X^1\cup X^2$ for $i\in\{1,2\}$,
and $\lambda_0: X^1\cap X^2\to X^1\cup X^2$ are inclusions.
\end{corollary}

\subsection{A spectral sequence of a morphism}

Let $f: X\rightarrow Y$ be a morphism of cubical sets.
The inverse images of singular cubes $\widetilde{y}$ of $Y$
form a diagram of cubical sets
$\{\overleftarrow{f}(y)\}_{\widetilde{y}\in \Box/Y}$ whose colimit is isomorphic to
 $X$. Applying the general theorem on the spectral sequence of a morphism
 \cite[Theorem 4.1]{X1991}, where we must
 take the category $\Ab^{op}$ instead of  $\mA$, we obtain the
 following statement with the help of Theorem \ref{main2}:

\begin{corollary}\label{spmor}
Let $f: X\rightarrow Y$ be a morphism of cubical sets and
$F$ a contravariant system on $X$. Then there is a first quadrant
spectral sequence
$$
E^2_{p,q} = \coLim_p^{{\Box}/Y}\{H_q(\overleftarrow{f}(y),
f_{y}^*(F))\}_{\widetilde{y}\in{{\Box}/Y}} \Rightarrow
H_{p+q}(X,F).
$$
\end{corollary}

\section{The conclusion}

We have established that the homology groups of a cubical set with coefficients
in a contravariant system can be considered as left derived functors of the colimit
with coefficients in this system.
The corollaries obtained show that the homology of cubical sets has many  properties
known earlier for simplicial, semisimplicial and precubical sets.

In the future, the applications of the main theorem for the study of homology
groups of cubical sets with coefficients in local systems, and also
for studying homology groups of directed topological spaces.

\end{document}